# A PROBABILISTIC ANALYSIS OF SOME TREE ALGORITHMS

By Hanène Mohamed and Philippe Robert

*INRIA*

In this paper a general class of tree algorithms is analyzed. It is shown that, by using an appropriate probabilistic representation of the quantities of interest, the asymptotic behavior of these algorithms can be obtained quite easily without resorting to the usual complex analysis techniques. This approach gives a unified probabilistic treatment of these questions. It simplifies and extends some of the results known in this domain.

**1. Introduction.** A splitting algorithm is a procedure that divides recursively into subsets an initial set of $n$ items until each of the subsets obtained has a cardinality strictly less than some fixed number $D$. These algorithms have a wide range of applications:

(a) Data structures. These are algorithms on data structures used to sort and search. They are sometimes referred to as divide and conquer algorithms. See [6] and [23] for a general presentation and [28, 39, 40] for their analysis with analytical methods.

(b) Communication networks. These algorithms are used to give a distributed access to a common communication channel that can transmit only one message per time unit. See [3, 10, 41].

(c) Distributed systems. Some algorithms use a splitting technique to select a subset of a set of identical communicating components. See [18, 36].

(d) Statistical tests. A test, performed on a set of individuals, indicates if at least one of these individuals has some characteristics (like a disease if this is blood testing). The purpose is to minimize the number of tests to identify individuals with the specified characteristic as quickly as possible. See [43].

Formally, a splitting algorithm can be described as follows:

SPLITTING ALGORITHM $\mathcal{S}(n)$
— TERMINATION CONDITION.
  If $n < D$   $\longrightarrow$ STOP.
— TREE STRUCTURE.
  If $n \geq D$, randomly divide $n$ into $n_1, \ldots, n_G$, with $n_1 + \cdots + n_G = n$ where $G$ is a random variable with some fixed distribution.
  $\longrightarrow$ APPLY $\mathcal{S}(n_1)$, $\mathcal{S}(n_2)$, ..., $\mathcal{S}(n_G)$.

1.1. *Description.* The algorithm starts with a set of $n$ items. This set is randomly split into $G$ subsets, the distribution of $G$ is given by $\mathbb{P}(G = \ell) = p_\ell$, where $(p_\ell)$ is a probability distribution on $\{2, 3, \ldots\}$. Now, conditionally on the event $\{G = \ell\}$, for $1 \leq i \leq \ell$, an item is sent into the $i$th subset with probability $V_{i,\ell}$, where $V_\ell = (V_{i,\ell}; 1 \leq i \leq \ell)$ is a random probability vector on $\{1, \ldots, \ell\}$. It can also be seen as a vector of random weights on the $\ell$ arcs of the branching procedure on which each of the $n$ items perform a random walk.

If $N_i$ is the cardinality of the $i$th subset, then, conditionally on the event $\{G = \ell\}$ and on the random variables $V_{1,\ell}, V_{2,\ell}, \ldots, V_{\ell,\ell}$, the distribution of the vector $(N_1, \ldots, N_\ell)$ is multinomial with parameter $n$ and $(V_{1,\ell}, V_{2,\ell}, \ldots, V_{\ell,\ell})$,

$$\mathbb{P}((N_1, \ldots, N_\ell) = (m_1, \ldots, m_\ell)) = \frac{n!}{m_1! m_2! \cdots m_\ell!} \prod_{k=1}^{\ell} (V_{k,\ell})^{m_k},$$

for $(m_i) \in \mathbb{N}^n$ such that $m_1 + \cdots + m_\ell = n$. If the $i$th subset, $1 \leq i \leq n$, is such that $N_i < D$, the algorithm stops for this subset. Otherwise, it is applied to the $i$th subset: a variable $G_i$, with the same distribution as $G$, is drawn and this $i$th subset is split into $G_i$ subsets, and so on.

Such a random splitting has been introduced by Devroye [7] where the asymptotic expansion of the depth of the associated tree is investigated.

EXAMPLES. (i) Knuth's algorithm. When $\mathbb{P}(G = 2) = 1$, $D = 2$ and $V_{1,2} \equiv V_{2,2} \equiv 1/2$, this is one of the oldest algorithms of this kind. It was analyzed by Knuth in 1973.

(ii) Symmetrical splitting algorithm. This is the case where $V_{i,n} \equiv 1/n$ for any $n \geq 2$ and $1 \leq i \leq n$.

(iii) $Q$-ary algorithm. If $\mathbb{P}(G = Q) = 1$ and $D = 2$, this is the $Q$-ary resolution algorithm with blocked arrivals analyzed by Mathys and Flajolet [30].

See also [7] for other examples. Quite naturally, such an algorithm can be graphically represented with a tree as shown by Figure 1.

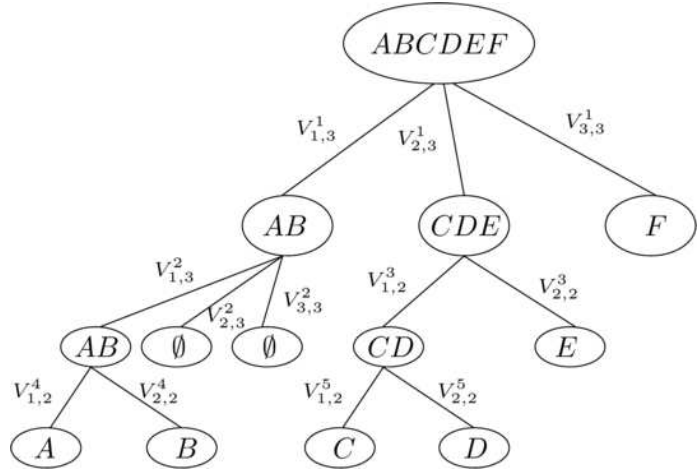

FIG. 1. *Splitting algorithm with $D=2$, two sets of random weights $(V_{1,2}, V_{2,2})$ and $(V_{1,3}, V_{2,3}, V_{1,3})$, $G$ a random variable with values in $\{2,3\}$ and the initial items $A$, $B$, $C$, $D$, $E$ and $F$.*

*Splitting measure.* As it will be seen in the following, the key characteristic of this splitting algorithm is a probability distribution $\mathcal{W}$ on $[0,1]$ defined with the branching distribution (the variable $G$) and the weights on each arc [the vector $(V_{1,G}, \ldots, V_{G,G})$]. The asymptotic behavior of the algorithm is expressed naturally in terms of the distribution $\mathcal{W}$.

DEFINITION 1. The *splitting measure* is the probability distribution $\mathcal{W}$ on $[0,1]$ defined by, for a nonnegative Borelian function $f$,

$$(1) \quad \int f(x)\mathcal{W}(dx) = \mathbb{E}\left(\sum_{i=1}^{G} V_{i,G} f(V_{i,G})\right) = \sum_{\ell=2}^{+\infty}\sum_{i=1}^{\ell} \mathbb{P}(G=\ell)\mathbb{E}(V_{i,\ell} f(V_{i,\ell})).$$

ASSUMPTION (A). Throughout the paper, it is assumed that, almost surely $G \geq 2$, and that there exists some $\delta > 0$ such that the relation

$$(A) \quad \sup_{\ell \geq 2} \sup_{1 \leq i \leq \ell} V_{i,\ell} \leq \delta < 1$$

holds almost surely, in particular $\mathcal{W}([0,\delta]) = 1$. These conditions imply in particular the nondegeneracy of the splitting mechanism.

DEFINITION 2. A splitting measure $\mathcal{W}$ is *exponentially arithmetic* if there exists some $\lambda > 0$ such that

$$\mathcal{W}(\{e^{-n\lambda} : n \geq 1\}) = 1,$$

and the largest $\lambda$ satisfying this relation is defined as the *exponential span* of $\mathcal{W}$.

If $A$ is some random variable with distribution $\mathcal{W}$, then $\mathcal{W}$ is exponentially arithmetic with exponential span $\lambda$ if and only if the distribution of $-\log(A)$ is arithmetic with span $\lambda$. See [14].

EXAMPLES. (i) Knuth's algorithm, $\mathbb{P}(G=2)=1$, $D=2$ and $V_{1,2} \equiv V_{2,2} \equiv 1/2$.

In this case

$$\mathcal{W}(dx) = \delta_{1/2},$$

where $\delta_x$ is the Dirac distribution at $x$ and $\mathcal{W}$ is exponentially arithmetic with exponential span $\log 2$.

(ii) Symmetrical splitting algorithm.

$$\mathcal{W}(dx) = \sum_{n\geq 2} P(G=n)\delta_{1/n},$$

the exponential span is $\log D$ where $D$ is the largest integer $p$ such that the support of the random variable $G$ is contained in $p\mathbb{N}$.

(iii) $Q$-ary algorithm, $\mathbb{P}(G=Q)=1$, $D=2$, $V_{1,Q}=p_1,\ldots,V_{Q,Q}=p_Q$:

$$\mathcal{W}(dx) = p_1\delta_{p_1} + p_2\delta_{p_2} + \cdots + p_Q\delta_{p_Q},$$

the distribution $\mathcal{W}$ is exponentially arithmetic if and only if all the real numbers $\log p_i/\log p_j$, $1\leq i<j\leq Q$, are rational.

*The cost of a splitting algorithm.* For such an algorithm, an important quantity is the number of operations required until the algorithm stops, that is, when all the subsets have a cardinality less than or equal to $D$. Denote by $R_n$ this quantity when the number of initial items is $n$; then clearly:

(a) $R_n = 1$ when $n < D$;
(b) for $n \geq D$,

$$R_n \stackrel{\text{dist.}}{=} 1 + R_{1,N_1^n} + \cdots + R_{G,N_G^n}, \quad (2)$$

where conditionally on the event $\{G=\ell\}$ and the random variables $V_{1,\ell}, V_{2,\ell}, \ldots, V_{\ell,\ell}$,

(1) the vector $(N_1^n,\ldots,N_\ell^n)$ has a multinomial distribution with parameter $n$ and $(V_{1,\ell}, V_{2,\ell},\ldots, V_{\ell,\ell})$;
(2) for $(p_i) \in \mathbb{N}^\ell$, the variables $R_{1,p_1},\ldots,R_{\ell,p_\ell}$ are independent;
(3) for $1\leq i\leq \ell$, the variable $R_{i,p_i}$ has the same distribution as $R_{p_i}$.

The variable $R_n$ is simply the number of nodes of the associated tree; see Figure 1.

1.2. *Unusual laws of large numbers.* Note that, since the splitting procedure is random, the variable $R_n$ is a *random variable*. With the language of communication networks, this quantity can be thought of as the total time to transmit $n$ initial messages. If $\mathbb{E}(R_n)$ is its expected value, $\mathbb{E}(R_n)/n$ is the average transmission time of one message among $n$. From a probabilistic point of view, it is natural to expect that the sequence $(R_n)$ satisfies a kind of law of large numbers, that is, that $(\mathbb{E}(R_n)/n)$ converges to some quantity $\alpha$. The constant $\alpha$ is, in some sense, the asymptotic average transmission time of a message. Curiously, this law of large numbers does not always hold. In some situations, the sequence $(\mathbb{E}(R_n)/n)$ does not converge at all and, moreover, exhibits an oscillating behavior.

When the splitting degree is constant and equal to 2 and $V_{1,2} \equiv V_{2,2} \equiv 1/2$ (the items are equally divided among the two subsets), these phenomena are quite well known. They have been analyzed using complex analysis techniques, functional transforms (and their associated inversion procedures) by Knuth [23], Flajolet, Gourdon and Dumas [13], Louchard and Prodinger [27] and many others. See [16, 28, 39] for a comprehensive treatment of this approach. See also [8] for a survey of the domain. Robert [38] proposed an alternative, elementary method to get the asymptotic behavior of some related oscillating sequences without using complex analysis.

When the splitting degree is constant and equal to $Q$ but the items are *not* equally divided among the subsets, studies are quite rare. Using complex analysis techniques, Fayolle, Flajolet and Hofri [12] obtained the asymptotic behavior of the associated sequence $(\mathbb{E}(R_n))$. Mathys and Flajolet [30] presented a sketch of a generalization of this study when $Q$ is arbitrary.

*Some alternative approaches.* (i) Some laws of large numbers have been proved by Devroye [8] in a quite general framework for various functionals of the associated trees. Talagrand's concentration inequalities are the main tools in this study. In our case, it would consist in proving that the distribution of the random variable $R_n/\mathbb{E}(R_n)$ is sharply (with an exponential decay) concentrated around 1. Results on limiting distributions such as central limit theorems do not seem to be accessible with this method.

(ii) Clement, Flajolet and Vallée [5] analyzed related algorithms in the more general context of dynamical systems. By using a Hilbertian setting, they showed that the first-order behavior of the algorithms is expressed in terms of the spectrum of a functional operator, the transfer operator. Getting explicit results in this way requires therefore a good knowledge of some eigenvalues of the transfer operator.

A dynamic version of this class of algorithms is investigated in [33]. The splitting procedure is the same but, in the language of branching processes, an immigration occurs at every leaf of the associated tree, that is, new items arrive every time unit. This dynamic feature complicates the problem. In this case, an additional probabilistic tool has to be used: an autoregressive process with moving average plays an important role.

### 1.3. Related problems.

*Fragmentation processes.* A continuous version of a splitting algorithm could be defined as follows: an initial mass of size $x$ is randomly split into several pieces and, at their turn, each of the pieces is randomly split.... A class of such models has been recently investigated. The fragmentation of each mass occurs after some independent exponential time with a parameter depending, possibly, on its mass. See [2, 32] and references therein. The problems considered are somewhat different: regularity properties of associated Markov processes, duality, rate of decay of individual masses, loss of mass, asymptotic distributions, and so on. A splitting algorithm is just a recursive fragmentation of an integer into integer pieces until each of the components has a size less than $D$. In a continuous setting, an analogue of the algorithms considered here would consist in stopping the fragmentation process of a mass as soon as its value is below some threshold $\varepsilon > 0$.

*Random recursive decompositions.* As it will be (easily) seen, a splitting algorithm can also be described as a random recursive splitting of the interval $[0,1]$. For example, in the case of a dyadic splitting, starting from the interval $[0,1]$, two subintervals $I_1$, $I_2$ are created and each of them is split at its turn and so on.

These random recursive decompositions have been considered from the point of view of the geometry of the boundary points by various authors, to express the Hausdorff dimension of this set of points in particular. Mauldin and Williams [31] and Waymire and Williams [42] considered decompositions of the interval $[0,1]$ which are not necessarily conservative, that is, when $|I_1| + |I_2| < 1$ holds with positive probability in the dyadic case.

Hambly and Lapidus [15] and Falconer [11] considered decompositions of the interval $[0,1]$ from the point of view of the lengths of the associated subintervals. The interval $[0,1]$ is represented by a nonincreasing sequence $(L_n)$ whose sum is 1. For $n \geq 1$, $L_n$ is the length of the $n$th largest interval of the decomposition. This description is similar to the classical representation of fragmentation processes. See [34].

In this setting, multiplicative cascades and martingales introduced by Mandelbrot [29] and Kahane and Peyrière [19] show up quite naturally. They have been analyzed quite extensively; see [1, 25] and references therein.

1.4. *An overview.* The purpose of this paper is twofold. First, it considers splitting algorithms with a random (and possibly unbounded) degree of splitting generalizing the previous studies in this domain. Second, and this is in fact the main point of the paper, it proposes a probabilistic approach that simplifies greatly the analysis of these algorithms. Moreover, as a by-product, a new direct representation of the asymptotic oscillating behavior is established.

The analysis proposed in this paper also starts from (2), but its treatment is significantly different from the analytic approach. After some transformation, (2) is interpreted as a probabilistic equation which is iterated by using appropriate independent random variables. Following the method of Robert [38], the next step is to perform a probabilistic de-Poissonnization and, by using Fubini's theorem conveniently, to represent the quantity $\mathbb{E}(R_n)$ by using a Poisson point process on the real line. The final, crucial step which differs from [38], consists in using the key renewal theorem to get the asymptotic behavior of the sequence $(\mathbb{E}(R_n))$.

The approach is elementary; its main advantage over the analytic treatment lies certainly in the use of the renewal theorem which gives *directly* the asymptotic behavior.

*Results of the paper.* Section 2 gives a useful representation for the average cost of the algorithm. The main result of the paper for the asymptotic cost is the following theorem in Section 3. This is a summary of Propositions 9 and 11.

THEOREM 3 (Asymptotics of the average cost). *For a splitting algorithm, under the condition*

$$\text{(3)} \qquad \int_0^1 \frac{|\log(y)|}{y} \mathcal{W}(dy) < +\infty,$$

— *if the splitting measure $\mathcal{W}$ is not exponentially arithmetic, then*

$$\text{(4)} \qquad \lim_{n \to +\infty} \frac{\mathbb{E}(R_n)}{n} = \frac{\mathbb{E}(G)}{(D-1) \int_0^1 |\log(y)| \mathcal{W}(dy)}.$$

— *If the splitting measure $\mathcal{W}$ is exponentially arithmetic with exponential span $\lambda > 0$, as $n$ gets large, the equivalence*

$$\text{(5)} \qquad \frac{\mathbb{E}(R_n)}{n} \sim F\left(\frac{\log n}{\lambda}\right)$$

*holds, where $F$ is the periodic function with period $1$ defined by, for $x \geq 0$,*

$$F(x) = \frac{\mathbb{E}(G)}{\int_0^1 |\log(y)| \mathcal{W}(dy)} \frac{\lambda}{1 - e^{-\lambda}}$$
$$\times \int_0^{+\infty} \exp\left(-\lambda\left\{x - \frac{\log y}{\lambda}\right\}\right) \frac{y^{D-2}}{(D-1)!} e^{-y} \, dy$$

*and $\{z\} = z - \lfloor z \rfloor$ is the fractional part of $z \in \mathbb{R}$.*

Condition (3) is not really restrictive since the variable $G$ is bounded in practice. This theorem covers and extends some of the results in this domain: for Knuth's algorithm [23] and for $Q$-ary algorithms with blocked arrivals [30], see Corollaries 10 and 12.

Furthermore, when there are asymptotic periodic oscillations, the periodic function $F$ involved is expressed directly and not in terms of its Fourier coefficients as is usually the case. The expression of $F$ generalizes the representation of [38] obtained for Knuth's algorithm.

The *distribution* of the sequence $(R_n)$ (and not only its average) is investigated in Section 4. For simplicity, only the case where the variable $G$ is constant and the variables $V_{\cdot,G}$ are equal to $1/G$ is considered. The purpose of this section is to show that the distribution of the Poisson transform of the sequence and, more generally, the distribution of most of the functionals of the associated tree, can be expressed quite simply in terms of Poisson processes and uniformly distributed random variables.

Two representations of the distribution of the Poisson transform as a functional of Poisson processes are derived. As a consequence, a law of large numbers is proved when the number of initial items has a Poisson distribution (Poisson transform). Moreover, the asymptotic oscillating behavior of the algorithm is proved as a consequence of a *standard* law of large numbers. These unusual laws of large numbers are, in the end, in the realm of classical laws of large numbers.

The central limit theorem is also proved with a similar method in this case. This is a classical result (see [28]); it is usually proved with complex analysis methods via quite technical estimations. It is proved here as a consequence of the *standard* central limit theorem for independent random variables. At the same time, a new representation of the asymptotic variance is obtained.

**2. General properties.** Throughout this paper, $(\mathcal{N}([0,x]))$ denotes a Poisson process with intensity 1; equivalently it can also be described as a nondecreasing sequence $(t_n)$ such that $(t_{n+1}-t_n)$ is a sequence of i.i.d. random variables exponentially distributed with parameter 1. For $x \geq 0$, the variable $\mathcal{N}([0,x])$ is simply the number of $t_n$'s in the interval $[0,x]$. See [20] for basic results on Poisson processes.

Equation (2) and the boundary conditions for the sequence $(R_n)$ are summarized in the following relation, for $n \geq 0$:

$$R_n \stackrel{\text{dist.}}{=} 1 + R_{1,N_1^n} + \cdots + R_{G,N_G^n} - G\mathbb{1}_{\{n<D\}},$$

therefore,

(6) $$R_n - 1 \stackrel{\text{dist.}}{=} \sum_{i=1}^{G}(R_{i,N_i^n}-1) + G\mathbb{1}_{\{n\geq D\}}.$$

DEFINITION 4. The *Poisson transform* of a nonnegative sequence $(a_n)$ is defined as

$$\sum_{n \geq 0} a_n \frac{x^n}{n!} e^{-x} = \mathbb{E}(a_{\mathcal{N}([0,x])}). \tag{7}$$

The following proposition gives useful representations of the Poisson transform of the sequence of $(\mathbb{E}(R_n))$.

PROPOSITION 5 [Poisson transform of the sequence $(R_n)$]. *For $x > 0$,*

$$\mathbb{E}(R_{\mathcal{N}([0,x])}) = 1 + \mathbb{E}(G) \mathbb{E}\left( \sum_{i=0}^{+\infty} \frac{1}{\prod_{k=1}^{i} W_k} \mathbb{1}_{\{t_D \leq x \prod_{k=1}^{i} W_k\}} \right), \tag{8}$$

*where $(W_i)$ is an i.i.d. sequence of random variables with distribution $\mathcal{W}$.*

PROOF. If $n$ is a Poisson random variable with parameter $x$, the splitting property of Poisson variables (see [20], e.g.) shows that, conditionally on the event $\{G = \ell\}$ and on the variables $V_{1,\ell}, \ldots, V_{\ell,\ell}$, the variables $N_i^n$, $1 \leq i \leq \ell$, are independent and $N_i^n$ has a Poisson distribution with parameter $xV_{i,l}$. Consequently, for $x > 0$, if

$$\Phi(x) \stackrel{\text{def.}}{=} \frac{\mathbb{E}(R_{\mathcal{N}([0,x])}) - 1}{x \mathbb{E}(G)}, \tag{9}$$

it is easily checked that $\mathbb{E}(G)\Phi(x) \to R_1 - R_0 = 0$ as $x \searrow 0$.

Since $\{\mathcal{N}([0,x]) \geq D\} = \{t_D \leq x\}$, (6) gives the relation

$$\Phi(x) = \sum_{\ell=2}^{+\infty} \mathbb{P}(G = \ell) \mathbb{E}\left( \sum_{i=1}^{\ell} V_{i,\ell} \Phi(xV_{i,\ell}) \right) + \frac{1}{x} \mathbb{P}(t_D \leq x). \tag{10}$$

Equation (10) can then be rewritten as

$$\Phi(x) = \mathbb{E}(\Phi(xW_1)) + \mathbb{E}\left( \frac{1}{x} \mathbb{1}_{\{t_D \leq x\}} \right). \tag{11}$$

The iteration of (11) shows that, for $n \geq 1$,

$$\Phi(x) = \mathbb{E}\left( \Phi\left( x \prod_{k=1}^{n} W_k \right) \right) + \mathbb{E}\left( \sum_{i=0}^{n-1} \frac{1}{x \prod_{k=1}^{i} W_k} \mathbb{1}_{\{t_D \leq x \prod_{k=1}^{i} W_k\}} \right).$$

The assumption on the variable $G$ and the sequence of vectors $(V_n)$ implies that, almost surely, the sequence $(\prod_{k=1}^{n} W_k)$ converges to 0. The function $\Phi$ can thus be represented as

$$\Phi(x) = \mathbb{E}\left( \sum_{i=0}^{+\infty} \frac{1}{x \prod_{k=1}^{i} W_k} \mathbb{1}_{\{t_D \leq x \prod_{k=1}^{i} W_k\}} \right).$$

The proposition has been proved. □

From now on, throughout the paper, $(W_i)$ will denote an i.i.d. sequence of random variables on $[0,1]$ with distribution $\mathcal{W}$.

PROPOSITION 6 (Probabilistic de-Poissonnization). *For $n \geq D$, then*

$$\mathbb{E}(R_n) = 1 + \mathbb{E}(G)\mathbb{E}\left(\sum_{i=0}^{T(U_{(n)}^D)-1} \frac{1}{\prod_{k=1}^i W_k}\right), \tag{12}$$

*where, for $0 < y < 1$,*

$$T(y) = \inf\left\{i \geq 1 : \prod_{k=1}^i W_k < y\right\}$$

*and $U_{(n)}^D$ is the $D$th smallest variable of $n$ independent, uniformly distributed random variables on $[0,1]$ independent of $(W_i)$.*

PROOF. For $x > 0$, by decomposing with respect to the number of points of the Poisson process $(\mathcal{N}(t))$ in the interval $[0,x]$, one gets, for $0 < \alpha \leq 1$,

$$\mathbb{P}(t_D \leq x\alpha) = \sum_{n=D}^{+\infty} \mathbb{P}(t_D \leq x\alpha, \mathcal{N}([0,x]) = n)$$

$$= \sum_{n=D}^{+\infty} \mathbb{P}(t_D \leq x\alpha | \mathcal{N}([0,x]) = n)\mathbb{P}(\mathcal{N}([0,x]) = n).$$

For $n \geq D$, conditionally on the event $\{\mathcal{N}([0,x]) = n\}$, the variable $t_D$ has the same distribution as the $D$ smallest random variable of $n$ uniformly distributed random variables on $[0,x]$. When $x = 1$, denote by $U_{(n)}^D$ a variable with this conditional distribution. Clearly, by homogeneity, the variable $(t_D | \mathcal{N}([0,x]) = n)$ has the same distribution as $xU_{(n)}^D$. Finally, one gets the identity

$$\mathbb{P}(t_D \leq x\alpha) = \sum_{n=D}^{+\infty} \mathbb{P}(U_{(n)}^D \leq \alpha)\frac{x^n}{n!}e^{-x}$$

$$= \mathbb{E}\left(\sum_{n=D}^{+\infty} \mathbb{1}_{\{U_{(n)}^D \leq \alpha\}}\frac{x^n}{n!}e^{-x}\right).$$

By using the independence of the sequence $(W_i)$ and $t_D$ in (8), the last identity gives the relation

$$\mathbb{E}(R_{\mathcal{N}([0,x])}) = 1 + \mathbb{E}(G)\mathbb{E}\left(\sum_{i=0}^{+\infty} \frac{1}{\prod_{k=1}^i W_k} \sum_{n=D}^{+\infty} \mathbb{1}_{\{U_{(n)}^D \leq \prod_{k=1}^i W_k\}}\frac{x^n}{n!}e^{-x}\right).$$

By Fubini's theorem and writing $1 = \exp(x)\exp(-x)$, this expression can be rewritten as

$$\sum_{n=0}^{D-1} \frac{x^n}{n!} e^{-x} + \sum_{n=D}^{+\infty} \left(1 + \mathbb{E}(G)\mathbb{E}\left(\sum_{i=0}^{+\infty} \frac{1}{\prod_{k=1}^{i} W_k} \mathbb{1}_{\{U_{(n)}^D \leq \prod_{k=1}^{i} W_k\}}\right)\right) \frac{x^n}{n!} e^{-x}.$$

The identification of (7) of $\mathbb{E}(R_{\mathcal{N}([0,x])})$ and the last identity gives (12). $\square$

COROLLARY 7 (Symmetrical $Q$-ary algorithm). When $\mathbb{P}(G = Q) = 1$ holds and $V_{i,Q} \equiv 1/Q$, for $i = 1, \ldots, Q$, then for $n \geq D$,

$$(13) \qquad \mathbb{E}(R_n) = 1 + \frac{Q}{Q-1}(\mathbb{E}(Q^{\lceil -\log_Q U_{(n)}^D \rceil}) - 1)$$

with, for $0 \leq x \leq 1$,

$$\mathbb{P}(U_{(n)}^D > x) = \sum_{k=0}^{D-1} \binom{n}{k} x^k (1-x)^{n-k}.$$

From (13), by using the fact that $nU_{(n)}^D$ converges in distribution as $n$ tends to infinity, it is not difficult to get the asymptotic behavior of $\mathbb{E}(R_n)$. The general case, (12), is slightly more complicated. One has to study the asymptotics of the series inside the expectation.

2.1. *A functional integral equation.* If $R(x) = \mathbb{E}(R_{\mathcal{N}(]0,x])})$ denotes the expected value of the Poisson transform of the sequence $(R_n)$, then (6) gives the relation

$$R(x) = 1 + \sum_{\ell=2}^{+\infty} \mathbb{P}(G = \ell)\mathbb{E}\left(\sum_{i=1}^{\ell} R(xV_{i,\ell})\right) - \mathbb{E}(G)\mathbb{P}(t_D \geq x);$$

by denoting

$$h(x) = 1 - \mathbb{E}(G)\int_x^{+\infty} \frac{u^{D-1}}{(D-1)!} du,$$

it is easy to see that the above identity can be written as the following integral equation:

$$(14) \qquad R(x) = \int_0^{+\infty} R(xu) \frac{\mathcal{W}(du)}{u} + h(x).$$

Recall that $\mathcal{W}$ is some probability distribution on the interval $[0, 1]$. For the $Q$-ary protocol considered by Mathys and Flajolet [30], this equation is

$$R(x) = \sum_{i=1}^{Q} R(xp_i) + h(x).$$

It is analyzed by considering the Mellin transform $R^*(s)$ of $R(x)$ on some vertical strip $\mathcal{S}$ of $\mathbb{C}$,

$$R^*(s) = \int_0^{+\infty} R(u) u^{s-1}\, du, \qquad s \in \mathcal{S},$$

which, in this case, is given by

$$R^*(s) = h^*(s) \bigg/ \left(1 - \sum_{i=1}^Q \frac{1}{p_i^s}\right).$$

The analytical approach consists in analyzing the poles of $R^*(s)$ on the right-hand side of $\mathcal{S}$, basically the solutions with positive real part of the equation

$$p_1^{-s} + p_2^{-s} + \cdots + p_Q^{-s} = 1.$$

Then, by inverting the Mellin transform and using complex analysis techniques, the asymptotic behavior of $(R(x))$ at infinity is described in terms of these poles. The final step, an analytic inversion of the Poisson transform together with technical estimates, establishes a relation between the asymptotic behaviors of the function $x \to R(x)$ and of the sequence $(R_n)$.

In the general case considered here, (14) gives the following expression for the Mellin transform of $(R(x))$:

$$R^*(s) = h^*(s) \bigg/ \left(1 - \int_0^{+\infty} \frac{1}{u^{s+1}} \mathcal{W}(du)\right).$$

An analogue of the analytic approach would start with the study of the roots $s \in \mathbb{C}$, $\Re(s) \geq 0$, of the equation

(15) $$\int_0^{+\infty} \frac{1}{u^{s+1}} \mathcal{W}(du) = 1,$$

and, if possible, proceed with successive inversions of Mellin transform and Poisson transform.

As it will be seen, our direct approach reduces to the minimum the technical apparatus required for such an analysis. The Poisson transform of $(R_n)$ is also used in our method, but it is conveniently represented [see (8)] so that it can be right away inverted to give an explicit expression (12) for $\mathbb{E}(R_n)$ which will give directly the asymptotic behavior of the sequence $(\mathbb{E}(R_n))$.

Interestingly, $(L_n)$ denotes the nonincreasing sequence of the lengths of the subintervals of $[0,1]$ associated to the splitting procedure (see Section 1.3). The *zeta function* of the string $(L_n)$ is defined as the meromorphic function

$$\zeta(s) = \sum_{n \geq 1} L_n^s, \qquad s \in \mathbb{C},$$

see [15] and [24]. It is not difficult to see that the relation

$$\mathbb{E}(\zeta(s)) = \int_0^{+\infty} u^s \mathcal{W}(du) \bigg/ \bigg(1 - \int_0^{+\infty} u^s \mathcal{W}(du)\bigg)$$

holds. In particular, the poles of the zeta function of the associated random recursive string can be expressed in terms of the solutions of (15).

## 3. Analysis of the asymptotic average cost.

*An associated random walk.* If $(W_i)$ is an i.i.d. sequence with common distribution $\mathcal{W}$ defined by (1), the sequence $(B_i) = (-\log(W_i))$ is an i.i.d. sequence of nonnegative random variables. The random walk $(S_n)$ is associated to $(B_i)$,

$$S_n = B_1 + B_2 + \cdots + B_n, \qquad n \geq 0.$$

As it will be seen, the asymptotic behavior of the splitting algorithm depends a great deal on the distribution of $(B_i)$. For $x > 0$, the crossing time $\nu_x$ of level $x$ by $(S_n)$ is defined as

$$\nu_x = \inf\{n : S_n > x\}.$$

For $0 < y < 1$, the variable $T(y)$ of Proposition 6 is simply $\nu_{-\log(y)}$. See [9, 14] for the main results concerning renewal theory used in the following.

If $\Psi$ is defined as

$$\Psi(x) = \mathbb{E}\bigg(\sum_{i=0}^{\nu_x - 1} \exp\bigg(\sum_{k=1}^{i} B_k\bigg)\bigg), \qquad x > 0,$$

then by (12),

(16) $$\mathbb{E}(R_n) = 1 + \mathbb{E}(G)\mathbb{E}[\Psi(-\log(U_{(n)}^D))].$$

It is clear that $-\log(U_{(n)}^D)$ converges in distribution to $+\infty$ as $n$ goes to infinity. The asymptotic behavior of $\Psi$ at infinity is first analyzed; this function can be rewritten as

(17) $$\Psi(x)e^{-x} = \mathbb{E}\bigg(\sum_{i=0}^{\nu_x - 1} e^{S_i - x}\bigg) = \mathbb{E}\bigg(\sum_{i=1}^{\nu_x} e^{S_{\nu_x - i} - x}\bigg).$$

3.1. *The nonarithmetical case.* In this part, it is assumed that the distribution of $W_1$ is not exponentially arithmetic. See Definition 2.

LEMMA 8. *Under the condition*

$$\mathbb{E}\bigg(\frac{|\log(W_1)|}{W_1}\bigg) = \int_0^1 \frac{|\log(x)|}{x} \mathcal{W}(dx) < +\infty,$$

*the relation*

$$\sup_{x \geq 0} \mathbb{E}(e^{S_{\nu_x} - x}) < +\infty$$

*holds.*

PROOF. Lorden's inequalities (see [4, 26]) show that, for any $p \geq 0$,

$$\sup_{x \geq 0} \mathbb{E}((S_{\nu_x} - x)^p) \leq \frac{p+2}{(p+1)\mathbb{E}(B_1)}\mathbb{E}(B_1^{p+1});$$

thus, one gets the relation

$$\sup_{x \geq 0} \mathbb{E}(e^{S_{\nu_x} - x}) \leq \frac{1}{\mathbb{E}(B_1)} \int_0^{+\infty} (u+2)e^u \mathbb{P}(B_1 \geq u)\,du$$

$$= \mathbb{E}((B_1 + 1)e^{B_1}) - 1$$

$$= \mathbb{E}\left(\frac{-\log(W_1) + 1}{W_1}\right) - 1 < +\infty. \qquad \square$$

For $i > 1$, the renewal theorem shows that, when $x$ goes to infinity, the variable $S_{\nu_x - i} - x$ converges in distribution to $-(\tau^* + \tau_1 + \tau_2 + \cdots + \tau_{i-1})$, where the variables $(\tau_n)$ are i.i.d. distributed as $B_1$ and independent of $\tau^*$ whose distribution is given by

$$\mathbb{E}(f(\tau^*)) = \frac{1}{\mathbb{E}(B_1)} \int_0^{+\infty} f(u)\mathbb{P}(B_1 \geq u)\,du,$$

for any nonnegative Borelian function on $\mathbb{R}$. By Assumption (A), the increments of the random walk $(S_n)$ are bounded below by $-\log(\delta)$, therefore one gets the relation, for $1 < K \leq \nu_x$,

$$\sum_{i=K}^{\nu_x} e^{S_{\nu_x - i} - x} = e^{S_{\nu_x} - x} \sum_{i=K}^{\nu_x} e^{S_{\nu_x - i} - S_{\nu_x}} \leq e^{S_{\nu_x} - x} \frac{\delta^K}{1 - \delta}.$$

From Lemma 8 and (17), one deduces then

$$\lim_{x \to +\infty} \Psi(x)e^{-x} = \mathbb{E}\left(\sum_{i=1}^{+\infty} \exp(-\tau^* - \tau_1 - \tau_2 - \cdots - \tau_{i-1})\right)$$

(18)
$$= \frac{1 - \mathbb{E}(\exp(-\tau_1))}{\mathbb{E}(\tau_1)} \times \frac{1}{1 - \mathbb{E}(\exp(-\tau_1))}$$

$$= \frac{1}{-\mathbb{E}(\log(W_1))},$$

since the density of $\tau^*$ on $\mathbb{R}_+$ is given by

$$\mathbb{P}(\tau_1 \geq x)/\mathbb{E}(\tau_1), \qquad x \geq 0.$$

PROBABILISTIC ANALYSIS OF TREE ALGORITHMS 15PROPOSITION 9 (Convergence of averages). *If the distribution of $W_1$ is not exponentially arithmetic and such that*

$$\mathbb{E}\left(\frac{|\log(W_1)|}{W_1}\right) < +\infty,$$

*then the following convergence holds:*

$$\lim_{n \to +\infty} \frac{\mathbb{E}(R_n)}{n} = \frac{\mathbb{E}(G)}{(D-1)\mathbb{E}(-\log W_1)}.$$

PROOF. Equation (16) gives that, for $n \geq 1$,

$$\frac{\mathbb{E}(R_n)}{n} = \frac{1}{n} + \mathbb{E}(G)\mathbb{E}\left(\Psi[-\log(U_{(n)}^D)]\exp(\log(U_{(n)}^D))\frac{1}{nU_{(n)}^D}\right).$$

As $n$ goes to infinity the variable $nU_{(n)}^D$ converges in distribution to a random variable $t_D$ which is a sum of $D$ i.i.d. exponential random variables with parameter 1; furthermore,

$$\lim_{n \to +\infty} \mathbb{E}\left(\frac{1}{nU_{(n)}^D}\right) = \mathbb{E}\left(\frac{1}{t_D}\right) = \frac{1}{D-1}.$$

For $\varepsilon > 0$, there exists $K$ such that, for $x > K$, $|\Psi(x)\exp(-x) + 1/\mathbb{E}(\log W_1)| < \varepsilon$; if $C$ denotes the supremum of $x \to \Psi(x)\exp(-x)$ on $\mathbb{R}_+$, then

$$\left|\mathbb{E}\left(\Psi[-\log(U_{(n)}^D)]\exp(\log(U_{(n)}^D))\frac{1}{nU_{(n)}^D}\right) - \frac{1}{(D-1)\mathbb{E}(-\log W_1)}\right|$$

(19)
$$\leq \varepsilon \mathbb{E}\left(\frac{1}{nU_{(n)}^D}\right) + \left(C + \frac{1}{\mathbb{E}(-\log W_1)}\right)\mathbb{E}\left(\mathbb{1}_{\{U_{(n)}^D > \exp(-K)\}}\frac{1}{nU_{(n)}^D}\right)$$

$$+ \frac{1}{\mathbb{E}(-\log W_1)}\left|\mathbb{E}\left(\frac{1}{nU_{(n)}^D}\right) - \frac{1}{D-1}\right|.$$

For $K_2 > 0$,

$$\limsup_{n \to +\infty} \mathbb{E}\left(\mathbb{1}_{\{U_{(n)}^D > \exp(-K)\}}\frac{1}{nU_{(n)}^D}\right) \leq \limsup_{n \to +\infty} \mathbb{E}\left(\mathbb{1}_{\{nU_{(n)}^D > K_2 \exp(-K)\}}\frac{1}{nU_{(n)}^D}\right)$$

$$= \mathbb{E}\left(\mathbb{1}_{\{t_D > K_2 \exp(-K)\}}\frac{1}{t_D}\right),$$

and this term goes to 0 as $K_2$ tends to infinity. One concludes that the right-hand side of (19) is arbitrarily small as $n$ goes to infinity. The proposition is proved. □

COROLLARY 10. (1) *Q-ary protocol with blocked arrivals.* When $D = 2$, $G \equiv Q$ and $V_{i,Q} = p_i$ for $1 \leq i \leq Q$, then, if at least one of the real numbers $\log p_i / \log p_1$, $2 \leq i \leq Q$, is not rational, the convergence

$$\lim_{n \to +\infty} \frac{\mathbb{E}(R_n)}{n} = \frac{Q}{\sum_{i=1}^{Q} -p_i \log p_i}$$

holds.

(2) *Symmetrical case.* If $G$ is not a degenerated random variable such that $\mathbb{E}(\log G) < +\infty$ and, for $\ell \geq 2$ and $1 \leq i \leq \ell$, $V_{i,\ell} = 1/\ell$, then

$$\lim_{n \to +\infty} \frac{\mathbb{E}(R_n)}{n} = \frac{\mathbb{E}(G)}{(D-1)\mathbb{E}(\log G)}.$$

3.2. *The arithmetical case.* It is assumed that the distribution of $W_1$ is exponentially arithmetic with exponential span $\lambda > 0$. The law of $-\log(W_1)/\lambda$ is a probability distribution on $\mathbb{N}$. For $i \geq 1$, one defines $C_i = B_i/\lambda = -\log(W_i)/\lambda$. In the arithmetic case, the integer-valued random walk associated to $(C_i)$ plays the key role, much in the same way as for $(S_n)$ in the nonarithmetic case. By denoting

$$\tau_n = \inf\left\{k \geq 1 : \sum_{i=1}^{k} C_i \geq n\right\},$$

(17) can be rewritten as, for $x \geq 0$,

$$\Psi(x) e^{-\lambda \lceil x/\lambda \rceil} = \mathbb{E}\left[\sum_{i=1}^{\tau_{\lceil x/\lambda \rceil}} \exp\left(\lambda \left(\sum_{k=1}^{\tau_{\lceil x/\lambda \rceil} - i} C_k - \lceil x/\lambda \rceil\right)\right)\right],$$

where $\lceil y \rceil = \inf\{n \in \mathbb{N} : n > y\}$ for $y \geq 0$. By using the discrete renewal theorem, for $i \geq 1$, as $n$ goes to infinity, the variable $C_1 + \cdots + C_{\tau_n - i} - n$ converges in distribution to $-(C_1^* + C_2 + \cdots + C_i)$, where $C_1^*$ is an independent random variable whose distribution is given by

$$\mathbb{P}(C_1^* = n) = \frac{1}{\mathbb{E}(C_1)} \mathbb{P}(C_1 \geq n), \qquad n \geq 1.$$

With the same method as in the nonarithmetic case, if the variable $|\log(W_1)|/W_1$ is integrable, then

$$\lim_{x \to +\infty} \Psi(x) e^{-\lambda \lceil x/\lambda \rceil} = \frac{1}{\mathbb{E}(|\log(W_1)|)} \frac{\lambda e^{-\lambda}}{1 - e^{-\lambda}}.$$

PROPOSITION 11 (Asymptotic periodic oscillations). *If the distribution of $W_1$ is exponentially arithmetic with exponential span $\lambda > 0$, and such that*

$$\mathbb{E}\left(\frac{|\log(W_1)|}{W_1}\right) < +\infty,$$

*then, as n gets large, the equivalence*

$$\frac{\mathbb{E}(R_n)}{n} \sim F\left(\frac{\log n}{\lambda}\right)$$

*holds, where $F$ is the periodic function with period 1 defined by, for $x \geq 0$,*

$$F(x) = \frac{\mathbb{E}(G)}{\mathbb{E}(|\log(W_1)|)} \frac{\lambda}{1-e^{-\lambda}}$$

$$\times \int_0^{+\infty} \exp\left(-\lambda\left\{x - \frac{\log y}{\lambda}\right\}\right) \frac{y^{D-2}}{(D-1)!} e^{-y} \, dy$$

*and $\{x\} = x - \lfloor x \rfloor$.*

PROOF. For $n \geq 1$, if $\lceil x \rceil = \lfloor x \rfloor + 1$,

$$\frac{1}{n}\mathbb{E}[\Psi(-\log(U_{(n)}^D))]$$

$$= \mathbb{E}\left[\Psi(-\log U_{(n)}^D) e^{-\lambda\lceil -\log(U_{(n)}^D/\lambda)\rceil} e^\lambda \exp\left(-\lambda\left\{\frac{-\log(U_{(n)}^D)}{\lambda}\right\}\right) \frac{1}{nU_{(n)}^D}\right];$$

since $nU_{(n)}^D$ converges in distribution to $t_D$ as $n$ goes to infinity, with the same method as in the proof of Proposition 9, one gets the equivalences

$$\frac{1}{n}\mathbb{E}[\Psi(-\log(U_{(n)}^D))] \times \mathbb{E}(|\log(W_1)|) \frac{1-e^{-\lambda}}{\lambda}$$

$$\sim E\left[\exp\left(-\lambda\left\{\frac{\log(n)}{\lambda} - \frac{\log(nU_{(n)}^D)}{\lambda}\right\}\right) \frac{1}{nU_{(n)}^D}\right]$$

$$= E\left[\exp\left(-\lambda\left\{\frac{\log(n)}{\lambda} - \frac{\log t_D}{\lambda}\right\}\right) \frac{1}{t_D}\right].$$

One concludes by using (16). □

COROLLARY 12 (*Q-ary protocol with blocked arrivals*). *When $D = 2$, $G \equiv Q$ and $V_{i,Q} = p_i$ for $1 \leq i \leq Q$, then, if all the real numbers $\log p_i / \log p_1$, $2 \leq i \leq Q$, are rational, the equivalence*

$$\frac{\mathbb{E}(R_n)}{n} \sim F\left(\frac{\log n}{\lambda}\right)$$

*holds, where $F$ is the periodic function with period 1 defined by, for $x \geq 0$,*

$$F(x) = \frac{Q}{-\sum_{i=1}^Q p_i \log p_i} \frac{\lambda}{1-e^{-\lambda}} \int_0^{+\infty} \exp\left(-\lambda\left\{x - \frac{\log y}{\lambda}\right\}\right) e^{-y} \, dy,$$

*where $\{x\} = x - \lfloor x \rfloor$ and $\lambda = \sup\{y > 0 : \forall i \in \{1, \ldots, Q\}, \log p_i \in y\mathbb{Z}\}$.*

**4. The distributions of the symmetrical $Q$-ary algorithm.** From now on, it is assumed that the branching degree of the splitting algorithm is constant, that is, $\mathbb{P}(G=Q)=1$, and uniform, $V_{i,Q} \equiv 1/Q$ for $1 \leq i \leq Q$. A set of $n \geq D$ items is randomly, equally divided into $Q$ subsets. From Proposition 11, it is known that

$$\mathbb{E}(R_n)/n \sim F_1(\log_Q n)$$

as $n$ goes to infinity, with

$$(20) \qquad F_1(x) = \frac{Q^2}{Q-1} \int_0^{+\infty} Q^{-\{x-\log_Q y\}} \frac{y^{D-2}}{(D-1)!} e^{-y}\, dy.$$

This is a typical case where a regular law of large numbers does not hold.

The purpose of this section is to strengthen the above convergence. The distribution of the Poisson transform of the sequence $(R_n)$, that is, the random variable $R_{\mathcal{N}(]0,x])}$, is investigated and not only its average as before. In particular it is shown that, for the Poisson transform, a *standard* law of large numbers can be used to prove the oscillating behavior of the algorithm. In other words, these uncommon laws of large numbers can be, in the end, expressed in a classical probabilistic setting.

NOTATION. Throughout the rest of the paper it is assumed that:

(1) $\mathcal{N}$ is a Poisson process with intensity 1 on $\mathbb{R}_+$. Another Poisson process will be used but in the two-dimensional space $[0,1] \times \mathbb{R}_+$.
(2) The variable $\mathcal{M}$ denotes a Poisson process on $[0,1] \times \mathbb{R}_+$ with intensity 1; this is a distribution of random points on $[0,1] \times \mathbb{R}_+$ with the following properties: if $\mathcal{M}(H)$ denotes the number of points that "fall" into the set $H \subset [0,1] \times \mathbb{R}_+$,

   (i) For $x \in [0,1] \times \mathbb{R}_+$, $M(\{x\}) \in \{0,1\}$.
   (ii) If $G$ and $H$ are disjoint subsets of $[0,1] \times \mathbb{R}_+$, the variables $\mathcal{M}(G)$ and $\mathcal{M}(H)$ are independent.
   (iii) The distribution of the variable $\mathcal{M}([a,b] \times [y,z])$ is Poisson with parameter $(b-a)(z-y)$ for $0 \leq a \leq b \leq 1$ and $0 \leq y \leq z$.

   Note that the random variables $\mathcal{N}([0,x])$ and $\mathcal{M}([0,1] \times [0,x])$ have a Poisson distribution with parameter $x$.
(3) The Poisson transform of the sequence $(R_n)$ is denoted by $\mathcal{R}(x)$, $x \geq 0$,

$$\mathcal{R}(x) \stackrel{\text{dist.}}{=} R_{\mathcal{N}([0,x])} \stackrel{\text{dist.}}{=} R_{\mathcal{M}([0,1] \times [0,x])}.$$

Its expectation is given by (8). This section is devoted to the study of the asymptotic behavior of the *distribution* of $\mathcal{R}(x)$.

4.1. *Laws of large numbers.* In this section it is proved that the Poisson transform of the sequence $(R_n)$ satisfies a strong law of large numbers. A nice representation of this transform as a functional of Poisson processes is first proved in the following proposition.

PROPOSITION 13. *The distribution of the Poisson transform $\mathcal{R}(x)$ of the sequence $(R_n)$ satisfies the following relations:*

$$(21) \quad \mathcal{R}(x) \stackrel{dist.}{=} \mathcal{R}_1(x) \stackrel{def.}{=} 1 + Q \sum_{p \geq 0} \sum_{k=0}^{Q^p-1} \phi_{\mathcal{N}}(xk/Q^p, x(k+1)/Q^p),$$

*where, for $0 \leq a \leq b$, $\phi_{\mathcal{N}}(a,b) = 1$ if $\mathcal{N}(]a,b]) \geq D$ and $0$ otherwise,*

$$(22) \quad \mathcal{R}(x) \stackrel{dist.}{=} \mathcal{R}_2(x) \stackrel{def.}{=} 1 + Q \sum_{p \geq 0} \sum_{k=0}^{Q^p-1} \mathbb{1}_{\{\mathcal{M}(]k/Q^p, (k+1)/Q^p] \times [0,x]) \geq D\}}.$$

Note that the function $x \to \mathcal{R}_2(x)$ is clearly nondecreasing. In particular, if $f$ is some nondecreasing function on $\mathbb{R}_+$, the same property holds for

$$x \to \mathbb{E}[f(\mathcal{R}(x))].$$

Representation (21) will be useful to get a strong law of large numbers on subsequences and also will be used to get the full convergence in distribution of $\mathcal{R}(x)/x$ as $x$ tends to infinity.

PROOF OF PROPOSITION 13. By the splitting property of Poisson random variables, the recurrence relation (2) for the sequence $(R_n)$ can be expressed as

$$R_{\mathcal{N}(]0,x])} \stackrel{dist.}{=} 1 + \sum_{i=1}^{Q} R_{i, \mathcal{N}(]x(i-1)/Q, xi/Q])} - Q \mathbb{1}_{\{\mathcal{N}(]0,x]) < D\}},$$

for $x \geq 0$. If, for $0 \leq a < b$,

$$(23) \quad \Phi(a,b) = \frac{1}{Q}(R_{\mathcal{N}(]a,b])} - 1),$$

the last equation can be rewritten as

$$(24) \quad \Phi(0,x) = \sum_{i=1}^{Q} \Phi_i\left(\frac{i-1}{Q}x, \frac{i}{Q}x\right) + \phi_{\mathcal{N}}(0,x),$$

with an obvious notation with the subscripts $i$ for $\Phi$. By iterating this relation, one gets that, almost surely, the expansion

$$(25) \quad \Phi(0,x) = \sum_{p \geq 0} \sum_{k=0}^{Q^p-1} \phi_{\mathcal{N}}\left(\frac{k}{Q^p}x, \frac{k+1}{Q^p}x\right)$$

holds. The function $\Phi(0,x)$ is just the sum of the function $\phi_\mathcal{N}$ on the $Q$-adic intervals of $[0,x]$.

Equation (22) is proved in the same way. □

*Representation of some of the functionals of the associated tree.* When $\mathcal{N}([0,x])$ items are at the root of the associated tree, the total number of nodes of the tree $R_{\mathcal{N}([0,x])}$ is not the only quantity that can be represented, by (21) in terms of the Poisson process $\mathcal{N}$.

The *maximal depth* $M(x)$ of the associated tree when there are $\mathcal{N}([0,x])$ items at the top of the tree can be expressed as a functional of the Poisson process

$$M(x) = \max\{p \geq 1 : \exists k, 0 \leq k < Q^{p-1} - 1, \mathcal{N}(]k/Q^{p-1}, (k+1)/Q^{p-1}]) \geq D\}.$$

The quantity

$$F(x) = \max\{p \geq 1 : \forall k, 0 \leq k < Q^{p-1} - 1, \mathcal{N}(]k/Q^{p-1}, (k+1)/Q^{p-1}]) \geq D\}$$

is the *number of full levels* of the tree. See [22]. Note that these quantities are directly related to classical occupancy problems. The *number of nodes at level* $p \geq 1$ is given by

$$Q \sum_{k=0}^{Q^{p-1}-1} \mathbb{1}_{\{\mathcal{N}(]k/Q^{p-1},(k+1)/Q^{p-1}]) \geq D\}}.$$

This is not, of course, an exhaustive list of the possible representations in terms of the Poisson process.

It is quite useful to think of splitting algorithms either in terms of trees or in terms of $Q$-adic subintervals of $[0,1]$. In a more general case, that is, when the splitting algorithm is not symmetrical, a representation similar to (21) can be obtained by using the associated random decomposition of the interval $[0,1]$ instead of the $Q$-adic decomposition. See [11].

*A strong law of large numbers.* Equation (24) shows that, if $N > 0$, the quantity $\Phi(0, yQ^N)$ is the sum of the $\Phi$ on the intervals $[yp, yp + y]$, $0 \leq p < Q^N$, and of $\phi_\mathcal{N}$ on the intervals $[ykQ^n, y(k+1)Q^n]$ contained in $[0, yQ^N]$, that is,

$$(26) \quad \Phi(0, yQ^N) = \sum_{p=0}^{Q^N-1} \Phi(yp, yp+y) + \sum_{n=1}^{N} \sum_{k=0}^{Q^{N-n}-1} \phi_\mathcal{N}(ykQ^n, y(k+1)Q^n).$$

By the independence properties of the Poisson process, the classical strong law of large numbers shows that, almost surely,

$$\lim_{N \to +\infty} \frac{1}{Q^N} \sum_{p=0}^{Q^N-1} \Phi(yp, yp+y) = \mathbb{E}(\Phi(0,y)) = \sum_{p \geq 0} Q^p \mathbb{E}(\phi_\mathcal{N}(0, y/Q^p))$$

$$= \sum_{p \geq 0} Q^p \mathbb{P}(\mathcal{N}(]0, y/Q^p]) \geq D),$$

by using (25), and for $n > 0$,

$$\lim_{N \to +\infty} \frac{1}{Q^N} \sum_{k=0}^{Q^{N-n}-1} \phi_{\mathcal{N}}(ykQ^n, y(k+1)Q^n) = \frac{1}{Q^n} \mathbb{E}(\phi_{\mathcal{N}}(0, yQ^n)).$$

Note that, for $0 < K < N$,

$$\sum_{n=K}^{N} \frac{1}{Q^N} \sum_{k=0}^{Q^{N-n}-1} \phi_{\mathcal{N}}(ykQ^n, y(k+1)Q^n) \leq \sum_{n=K}^{N} \frac{1}{Q^N} Q^{N-n} \leq \frac{1}{Q^{K-1}}.$$

The three last identities and decomposition (26) give that, almost surely,

$$\lim_{N \to +\infty} \frac{1}{yQ^N} \Phi(0, yQ^N) = \sum_{n \in \mathbb{Z}} \frac{1}{yQ^n} \mathbb{P}(\mathcal{N}(]0, yQ^n]) \geq D).$$

PROPOSITION 14 (Strong law of large numbers). *With the same notation as in Proposition* 13, *for* $0 < y \leq Q$, *almost surely,*

$$\lim_{N \to +\infty} \frac{\mathcal{R}_1(yQ^N)}{yQ^N} = Q \sum_{n \in \mathbb{Z}} \frac{1}{yQ^n} \mathbb{P}(\mathcal{N}(]0, yQ^n]) \geq D)$$

(27)
$$= Q \sum_{n \in \mathbb{Z}} \frac{1}{yQ^n} \int_0^{yQ^n} \frac{u^{D-1}}{(D-1)!} e^{-u} \, du$$

$$= F_1(\log_Q y),$$

*where* $F_1$ *is the periodic function defined by* (20).

As a by-product, the proposition establishes the intuitive (and classical) fact that the sequence $(\mathbb{E}(R_n)/n)$ and the function $x \to \mathbb{E}(R_{\mathcal{N}([0,x])})/x$ have the same asymptotic behavior at infinity. Note that if $G(y)$ is defined as the second term of (27), then the function $x \to G(Q^x)$ is clearly periodic with period 1.

PROOF OF PROPOSITION 14. Clearly, only the relation $F_1(\log_Q y) = G(y)$ has to be proved. For $n \in \mathbb{Z}$, if $t_D$ is the $D$th point of the Poisson process $\mathcal{N}$, then

$$\mathbb{P}(\mathcal{N}(]0, yQ^n]) \geq D) = \mathbb{P}(t_D \leq yQ^n) = \int_0^{+\infty} \mathbb{1}_{\{u \leq yQ^n\}} \frac{u^{D-1}}{(D-1)!} e^{-u} \, du.$$

By summing up these terms, with Fubini's theorem one gets

$$G(y) = Q\int_0^{+\infty} \sum_{n\in\mathbb{Z}} \frac{1}{yQ^n}\mathbb{1}_{\{u\leq yQ^n\}}\frac{u^{D-1}}{(D-1)!}e^{-u}\,du$$

$$= \frac{Q^2}{Q-1}\int_0^{+\infty} \frac{1}{yQ^{\lceil\log_Q(u/y)\rceil}}\frac{u^{D-1}}{(D-1)!}e^{-u}\,du$$

$$= \frac{Q^2}{Q-1}\int_0^{+\infty} \frac{1}{Q^{-\{\log_Q(u/y)\}}}\frac{u^{D-2}}{(D-1)!}e^{-u}\,du$$

$$= F_1(\log_Q y).$$

The proposition is proved. □

The following proposition establishes a weak law of large numbers for the Poisson transform of the sequence $(R_n)$. Devroye [8] obtained related results in a more general framework by using Talagrand's concentration inequalities.

THEOREM 15 (Law of large numbers). *The following convergence in distribution holds, for any $\varepsilon > 0$:*

$$\lim_{x\to+\infty}\mathbb{P}\left(\left|\frac{\mathcal{R}(x)}{xF_1(\log_Q x)} - 1\right| \geq \varepsilon\right) = 0,$$

*where $F_1$ is the function defined by (20).*

PROOF. For $x > 0$, one defines $N_x = \lfloor\log_Q x\rfloor$, $u_x = x/Q^{N_x}$ and, for $p \geq 1$, $z_x = \lfloor u_x p\rfloor Q^{N_x}/p$. Note that $\sup_{x\geq 1}|x/z_x - 1|$ converges to 0 as $p$ tends to infinity, hence by continuity of $F_1$,

$$\lim_{p\to+\infty}\sup_{x\geq 1}\left|\frac{xF_1(\log_Q x)}{z_x F_1(\log_Q z_x)} - 1\right| = 0.$$

Proposition 14 shows that for $p \geq 1$, almost surely, for $k$, $0 \leq k \leq p$,

$$\lim_{N\to+\infty}\frac{\mathcal{R}_1(y_k Q^N)}{y_k Q^N F_1(\log_Q y_k)} = 1,$$

with $y_k = k/p$. Therefore, if $p \geq 1$ is fixed, almost surely,

$$\lim_{x\to+\infty}\frac{\mathcal{R}_1(z_x)}{z_x F_1(\log_Q z_x)} = 1.$$

The monotonicity of the function $x \to \mathcal{R}_2(x)$ gives the relation

$$\mathcal{R}_2\left(\frac{\lfloor u_x p\rfloor}{p}Q^{N_x}\right) \leq \mathcal{R}_2(x).$$

One finally gets

$$\mathbb{P}\Big(\frac{\mathcal{R}(x)}{xF_1(\log_Q x)}<1-\varepsilon\Big)=\mathbb{P}\Big(\frac{\mathcal{R}_2(x)}{xF_1(\log_Q x)}<1-\varepsilon\Big)$$
$$\leq \mathbb{P}\Big(\frac{\mathcal{R}_2(z_x)}{xF_1(\log_Q x)}<1-\varepsilon\Big)$$
$$=\mathbb{P}\Big(\frac{\mathcal{R}_1(z_x)}{z_xF_1(\log_Q z_x)}\frac{z_xF_1(\log_Q z_x)}{xF_1(\log_Q x)}<1-\varepsilon\Big),$$

therefore,

$$\lim_{x\to+\infty}\mathbb{P}\Big(\frac{\mathcal{R}(x)}{xF_1(\log_Q x)}<1-\varepsilon\Big)=0.$$

The analogous inequality is obtained in the same way. The theorem is proved. □

4.2. *Central limit theorems.* For $N\geq 1$ and $0<x\leq Q$, with $\Phi$ defined by (23), the variance of the variable $\Phi(0,x)$ is first analyzed. The expansion (25) gives

$$[\Phi(0,x)-\mathbb{E}(\Phi(0,x))]^2=\sum_{p\geq 0}\sum_{k=0}^{Q^p-1}\sum_{p'\geq 0}\sum_{k'=0}^{Q^{p'}-1}\Delta_{k,p}(x)\Delta_{k',p'}(x),$$

with

$$\Delta_{k,p}(x)=\phi_\mathcal{N}\Big(\frac{k}{Q^p}x,\frac{k+1}{Q^p}x\Big)-\mathbb{E}\Big(\phi_\mathcal{N}\Big(0,\frac{1}{Q^p}x\Big)\Big).$$

The expected value of the variable $\Delta_{k,p}(x)\Delta_{k',p'}(x)$ is nonzero only if $p\leq p'$ and $kQ^{p'-p}\leq k'\leq (k+1)Q^{p'-p}-1$ or the symmetrical condition by exchanging $(p,k)$ and $(p',k')$:

$$\mathbb{E}[(\Phi(0,x)-\mathbb{E}(\Phi(0,x)))^2]$$
$$=\sum_{p\geq 0}\sum_{k=0}^{Q^p-1}\mathbb{E}[\Delta_{k,p}(x)^2]$$
$$+2\sum_{p\geq 0}\sum_{k=0}^{Q^p-1}\sum_{p'>p}\sum_{k'=kQ^{p'-p}}^{(k+1)Q^{p'-p}-1}\mathbb{E}[\Delta_{k,p}(x)\Delta_{k',p'}(x)].$$

By using the elementary identities

$$\mathbb{E}[\Delta_{k,p}(x)^2]=\mathbb{E}[\phi_\mathcal{N}(0,x/Q^p)](1-\mathbb{E}[\phi_\mathcal{N}(0,x/Q^p)])$$
$$=\mathbb{P}(t_D\leq x/Q^p)\mathbb{P}(s_D\geq x/Q^p),$$

$$\mathbb{E}[\Delta_{k,p}(x)\Delta_{k',p'}(x)] = \mathbb{E}[\phi_{\mathcal{N}}(0, x/Q^{p'})](1 - \mathbb{E}[\phi_{\mathcal{N}}(0, x/Q^p)])$$
$$= \mathbb{P}(t_D \leq x/Q^{p'})\mathbb{P}(s_D \geq x/Q^p),$$

where $t_D$ and $s_D$ are independent random variables with the same distribution as the $D$th point of the Poisson process $\mathcal{N}$, one gets the relation

$$\mathbb{E}[(\Phi(0,x) - \mathbb{E}(\Phi(0,x)))^2] = \sum_{p \geq 0} Q^p \mathbb{P}(t_D \leq x/Q^p \leq s_D)$$
$$+ 2 \sum_{p' > p \geq 0} Q^{p'} \mathbb{P}(t_D \leq x/Q^{p'}, x/Q^p \leq s_D).$$

By switching again the series and the expected values, one finally obtains

(28)
$$(Q-1)\mathbb{E}[(\Phi(0,x) - \mathbb{E}(\Phi(0,x)))^2]$$
$$= Q\mathbb{E}((Q^{\lfloor \log_Q(x/t_D) \rfloor} - Q^{\lfloor \log_Q(x/s_D) \rfloor})\mathbb{1}_{\{\lfloor \log_Q(x/t_D) \rfloor > \lfloor \log_Q(x/s_D) \rfloor\}})$$
$$+ 2Q\mathbb{E}((\lfloor \log_Q(x/t_D) \rfloor - \lfloor \log_Q(x/s_D) \rfloor - 1)^+ Q^{\lfloor \log_Q(x/t_D) \rfloor})$$
$$- 2\frac{Q}{Q-1}\mathbb{E}((Q^{\lfloor \log_Q(x/t_D) \rfloor} - Q^{\lfloor \log_Q(x/s_D) \rfloor + 1})^+),$$

where $a^+ = \max(a,0)$ for $a \in \mathbb{R}$. This identity gives the following proposition. A similar proposition has been proved by Jacquet and Régnier [17] and Régnier and Jacquet [37] in the case where $Q = D = 2$ but without symmetry conditions as is the case here. See also [28], Chapter 5.

PROPOSITION 16 (Asymptotic variance). *The variance of the Poisson transform of the sequence $(R_n)$ satisfies the following equivalence, as $x$ goes to infinity:*

$$\frac{1}{x} \text{Var}(\mathcal{R}(x)) \sim F_2(\log_Q(x)),$$

*where $F_2$ is the continuous periodic function with period 1 defined by, for $y \geq 0$,*

(29)
$$F_2(y) = \int_{\mathbb{R}_+^2} f_2(\{y - \log_Q(u)\}, \{y - \log_Q(v)\}, u, v)$$
$$\times \frac{u^{D-1}}{(D-1)!} \frac{v^{D-1}}{(D-1)!} e^{-(u+v)} \, du \, dv,$$

*with $\{z\} = z - \lfloor z \rfloor$ for $z \in \mathbb{R}$ and for $u > 0$, $v > 0$ and $y \in \mathbb{R}$,*

$$f_2(a,b,u,v) = \frac{Q}{Q-1}\left(\frac{Q^a}{u} - \frac{Q^b}{v}\right)\mathbb{1}_{\{\log_Q(v/u) + b > a\}}$$

$$+ \frac{2Q}{Q-1}(\log_Q(v/u) - a + b - 1)^+$$

$$- \frac{2Q}{(Q-1)^2}\left(\frac{Q^a}{u} - \frac{Q^{b+1}}{v}\right)^+$$

*where $z^+ = \max(z, 0)$.*

Note that a more detailed expansion of the variance could be obtained with (28).

PROPOSITION 17 (Central limit theorem for Poisson transform). *For $0 < y < Q$, as $N$ tends to infinity, the variable*

$$\frac{1}{\sqrt{Q^N}}(\mathcal{R}(yQ^N) - \mathbb{E}(\mathcal{R}(yQ^N)))$$

*converges in distribution to a Gaussian centered random variable with variance $yF_2(\log_Q y)$, where $F_2$ is defined by (29).*

PROOF. It is enough to prove the proposition for the variable $\Phi(0, x)$ defined by

$$\Phi(0, x) = \frac{1}{Q}(R_{\mathcal{N}(]0,x])} - 1).$$

Equation (26) gives, for $K \geq 1$,

$$\Phi(0, yQ^N) - \mathbb{E}(\Phi(0, yQ^N))$$

$$= \sum_{p=0}^{Q^N-1} [\Phi(yp, yp + y) - \mathbb{E}(\Phi(0, y))]$$

$$+ \sum_{n=1}^{K} \sum_{k=0}^{Q^{N-n}-1} [\phi_{\mathcal{N}}(ykQ^n, y(k+1)Q^n) - \mathbb{E}(\phi_{\mathcal{N}}(0, yQ^n))] + \Delta_K,$$

where $\Delta_K$ is the residual term of the series. By using the method to compute the variance, it is not difficult to establish that, for any $\varepsilon > 0$ there exists some $K > 0$ such that the expected value of $(\Delta_K/Q^N)^2$ is less than $\varepsilon$, for $N$ sufficiently large.

By regrouping the terms of the above equation according to the $Q$-adic intervals $[yk/Q^K, y(k+1)/Q^K]$ for $0 \leq k < Q^{N-K}$ and by using the independence properties of the Poisson process $\mathcal{N}$, the quantity $\Phi(0, yQ^N) - \mathbb{E}(\Phi(0, yQ^N)) - \Delta_K$ can be written as a sum of $Q^{N-K}$ independent identically distributed random variables. Therefore, the *classical* central limit theorem can be applied. The proposition is proved. □

4.3. *The distribution of the sequence* $(R_n)$. The following proposition describes the distribution of the variable $R_n$ in terms of $n$ i.i.d. uniformly distributed random variables on the interval $[0,1]$. This characterization is generally implicitly used to get various asymptotics describing the depth of the associated tree. See [28] and [35].

PROPOSITION 18. *For $n \geq 0$, the random variable $R_n$ has the same distribution as*

(30) $$R_n \stackrel{dist.}{=} 1 + Q \sum_{p \geq 0} \sum_{k=0}^{Q^p - 1} \phi_{\mathcal{U}_n}(]k/Q^p, (k+1)/Q^p]),$$

*where, for $0 \leq a \leq b \leq 1$, $\phi_{\mathcal{U}_n}(]a,b]) = 1$ if $\mathcal{U}_n(]a,b]) \geq D$ and $0$ otherwise. The variable $\mathcal{U}_n$ is the point measure on $[0,1]$ defined by*

$$\mathcal{U}_n = \delta_{U_1} + \delta_{U_2} + \cdots + \delta_{U_n},$$

$(U_1, \ldots, U_n)$ *are i.i.d. random variables uniformly distributed on $[0,1]$, in particular, $\mathcal{U}_n(]a,b])$ is the number of $U_i$'s in the interval $]a,b]$.*

PROOF. Assume that $\mathcal{N}$ is a Poisson process with parameter 1; by definition

$$(R_{\mathcal{N}(]0,x])} | \mathcal{N}(]0,x]) = n) \stackrel{dist.}{=} R_n.$$

Due to Proposition 13, the distribution of the Poisson transform $R_{\mathcal{N}(]0,x])}$ is expressed as a functional of the points of the Poisson process on the interval $[0,x]$. But, as in the proof of Proposition 6, conditionally on the event $\{\mathcal{N}(]0,x]) = n\}$, these points can be expressed as $xU_i$, $1 \leq i \leq n$, where $(U_i)$ are i.i.d. uniformly distributed random variables on the interval $[0,1]$. Equation (30) is thus a direct consequence of (21). □

INRIA
DOMAINE DE VOLUCEAU
B.P. 105
78153 LE CHESNAY CEDEX
FRANCE
E-MAIL: Hanene.Mohamed@inria.fr
        Philippe.Robert@inria.fr
URL: www-rocq.inria.fr/~robert